\documentclass{article}

\usepackage{amssymb}
\usepackage{graphicx}
\usepackage[all]{xy}

\usepackage[english,german,french,latin]{babel}

\def \m {{\medskip}}

\begin{document}

\begin{center}
{\bf{{\large{The Noether Theorems in Context}}}}

\m
Yvette Kosmann-Schwarzbach

\end{center}

\medskip

\noindent {\small {\it \selectlanguage{latin} \hfill {Methodi in hoc libro
    tradit{\ae}, non solum maximum 
    
\hfill esse usum in ipsa analysi, sed etiam eam 
 ad resolutionem 

\hfill problematum physicorum 
 amplissimum subsidium afferre.}}}

\hfill {\small Leonhard Euler [1744]}
\medskip

\medskip

``The methods described in this book are not
only of great use in analysis, but are also most helpful for the solution
of problems in
physics.'' 
Replacing  `in this book' by `in this article', the sentence that Euler wrote in the introduction to the first supplement of his treatise on the calculus of variations in 1744 applies equally well to Emmy Noether's 
\selectlanguage{german}``Invariante Variationsprobleme'', published in the {\it Nachrichten
von der K\"oniglichen Gesellschaft der Wissenschaften zu G\"ottingen,
Mathematisch-physikalische Klasse}\selectlanguage{english} in 1918.
\subsection*{Introduction}
In this talk,\footnote{This text is a revised version of the lecture I delivered at the international conference, ``The Philosophy and Physics of Noether's Theorems'', a centenary conference on the 1918 work of Emmy Noether, 
London, 5~October, 2018. To be published by Cambridge University Press (Nicholas Teh, James Read and Bryan Roberts, eds.)} I propose to sketch the contents of Noether's 1918 article,
``Invariante Variationsprobleme'', as it may be seen against the background
of the work of her predecessors and in the context of the debate on the conservation of energy that had arisen
in the general theory of 
relativity.\footnote{An English translation of Noether's article together with an account of her work and the history of its reception, from Einstein to Deligne, may be found in my book, {\it The Noether Theorems, Invariance and Conservation Laws in the Twentieth Century}, translated by Bertram~E. Schwarzbach [2010]. This book contains an extensive bibliography, only a small part of which is reproduced in the list of references below. It is an expanded version of my earlier book in French, {\it Les Th\'eor\`emes de Noether, Invariance et lois de conservation au XX$^e$ si\`ecle} [2004].}

Situating Noether's theorems on the invariant variational problems in their context requires a brief outline of the work of her
predecessors, and a description of her career, first in Erlangen, then in G\"ottingen.
Her 1918 article will be briefly summarised. I have endeavored to convey its contents in Noether's own vocabulary and notation with minimal recourse to more recent terminology.  Then I shall address these questions: how original was Noether's ``Invariante Variationsprobleme''? how modern were her use of Lie groups and her introduction of generalized vector fields? and
how influential was her article? To this end, I~shall sketch its reception from 1918 to 1970.
For many years, there was practically no recognition of either of these theorems. Then multiple references to ``the Noether theorem'' or ``Noether's theorem'' -- in the singular -- began to appear, either referring to her first theorem, in the publications of those mathematicians and mathematical physicists who were writing on mechanics -- who ignored her second theorem --, or to her second theorem by those writing on general relativity and, later, on gauge theory. 
I shall outline the curious transmission of her results, the history of the mathematical developments of her theory, and the ultimate recognition of the wide
applicability of ``the Noether theorems''.
To conclude, in the hope of dispelling various misconceptions, I~shall underline what Noether was {\it not}, and I shall reflect on the fortune of her theorems.

\subsection*{A family of mathematicians}
Emmy Noether was born to a Jewish family in 
Erlangen (Bavaria, Germany) in 1882. In a manuscript {\it curriculum vitae}, written for official purposes {\it circa} 1917, she described herself as 
``{of Bavarian nationality and Israelite
confession''.\footnote{Declaring one's religion was compulsory in Germany at the time.}
She died in Bryn Mawr (Pennsylvania) in the United States in 1935, after undergoing an operation. Why she had to leave Germany in 1933 to take up residence in America is clear from the chronology of the rise of the nazi regime in Germany and its access to power and has, of course, been told in the many accounts of her life that have been published,\footnote{The now classical biographies of Noether can be found in the book written by Auguste Dick [1970], translated into English in 1981, and in the volumes of essays edited by James W. Brewer and Martha K. Smith [1981], and by Bhama Srinivasan and Judith D. Sally [1983].} while numerous and sometimes fanciful comments have appeared in print and in the electronic media in recent years.

She was the daughter of the renowned mathematician, 
Max Noether (1844--1921), professor at the University of Erlangen. He had been a privatdozent, then an ``extraordinary professor'' in Heidelberg before moving to Erlangen in 1875, where he was eventually named an ``ordinary professor'' in 1888. 
Her brother, Fritz, was born in 1884, and studied mathematics and physics in Erlangen and Munich. He completed his doctorate in 1909 and became assistant to the professor of theoretical mechanics at the Technische Hochschule in Karlsruhe, where he submitted his {\it{Habilitation}} thesis in 1911. 
In 1922 he became professor in Breslau, from where he, too, was forced to leave in 1933. He emigrated to the Soviet Union and was appointed professor at the university of Tomsk.
Accused of being a German spy, he was jailed and shot in 1941.

\subsection*{The young Emmy Noether}
Emmy Noether first studied languages in order to become a teacher of French and English, a suitable profession for a young woman. But from 1900 on, she  studied mathematics, 
first in Erlangen, with her father, 
then audited lectures at the university.
For the winter semester in 1903--1904, she travelled to G\"ottingen to audit courses at that university. At that time, new regulations were introduced which enabled women to matriculate and take examinations. She then chose to enroll at the university of Erlangen, where she listed mathematics as her only course of study\footnote{On this, as well as on other oft repeated facts of Noether's biography, see Dick [1970], English translation, p. 14.},
and in 1907,
she completed her doctorate under the direction of 
{Paul Gordan} (1837--1912), a colleague of her father. Here I open a parenthesis: One should not confuse
the mathematician Paul Gordan, her 
\selectlanguage{german}\glqq Doktorvater\grqq\selectlanguage{english}, with the  physicist, Walter Gordon (1893--1939).
The ``Clebsch--Gordan coefficients'' in quantum mechanics bear the name of Noether's thesis adviser together with that of the physicist and mathematician 
Alfred Clebsch (1833--1872). However, the ``Klein--Gordon equation'' is named after Walter Gordon and the physicist Oskar Klein (1894--1977) who, in turn, should not be confused with the mathematician Felix~Klein about whom more will be said shortly.


\subsection*{Noether's 1907 thesis on invariant theory}

Noether's thesis at Erlangen University, entitled \selectlanguage{german}{``\"Uber die Bildung des Formensystems der tern\"aren
  biquadratischen Form''} \selectlanguage{english}(On the construction of the system of forms
of a ternary biquadratic form), dealt with the search for the invariants of those forms ({\it i.e.}, homogeneous
 polynomials) which are ternary ({\it{i.e.}}, in 3 variables) and biquadratic ({\it{i.e.}}, of degree 4).
An extract of her thesis appeared in 
he {\it Sitzungsberichte der
Physikalisch-medizinischen Societ\"at zu Erlangen} in 1907, and the complete text was 
published in 1908, in 
the \selectlanguage{german}{\emph{Journal 
f\"ur die reine und
  angewandte Mathematik}}\selectlanguage{english} (``Crelle's Journal").
She later distanced herself from her early work as having employed 
a needlessly computational approach to the problem.

After 1911, her work in algebra was influenced  by Ernst Fischer (1875--1954) who was appointed professor in Erlangen upon Gordan's retirement in 1910.
Noether's expertise in invariant theory 
revealed itself in publications in 1910, 1913, and 1915 that followed her thesis, and was later confirmed in the four articles on the invariants of finite groups that she published 
in 1916 in the {\it Mathematische Annalen}. 
She studied in particular 
the determination of bases of invariants that furnish an expansion with integral or rational coefficients of each invariant of the group, expressed as a linear combination of the invariants in the basis.

At Erlangen University from 1913 on, Noether occasionally substituted for her ageing
father, thus beginning to teach at the university level, but not under her own name. 

\subsection*{Noether's achievements}
Her achievement of 1918, whose centenary was duly celebrated in conferences in London and Paris, eventually became  a central result in both mechanics and field theory, and, more generally, in mathematical physics, though her role was rarely acknowledged before 1950 and, even then, it was only a truncated part of    her article that was cited. On the other hand, her articles on the theory of ideals and the representation theory of algebras published in the 1920's  made her world famous. Her role in the development of modern algebra was duly recognized by the mathematicians of the twentieth century, while they either considered her work on invariance principles to be an outlying and negligible part of her work or, more often, ignored it altogether.
In fact, the few early biographies of Noether barely mention her work on invariant variational problems, but both past and recent publications treat her fundamental contributions to modern algebra. 
I~shall not deal with them here. They are, and will no doubt continue to be celebrated by all mathematicians.

\subsection*{In G\"ottingen: Klein, Hilbert, Noether, and Einstein} 
In 1915, the great mathematicians, Felix Klein (1849--1925) and David Hilbert (1862--1943), invited Noether to G\"ottingen in the hope that her expertise in invariant theory would help them understand
some of the implications of
Einstein's 
newly formulated general theory of relativity.
In G\"ottingen, Noether took an active part in Klein's seminar. 
It was in her 1918 article that she solved a problem arising in the  general theory of relativity and proved ``the Noether theorems''. In particular, she proved and vastly generalized a conjecture made by Hilbert
concerning the nature of 
the law of conservation of energy.
Shortly afterwards, she returned to pure algebra.

At the invitation of Hilbert, Einstein had come to G\"ottingen in early July 1915 to deliver a series of lectures on the general theory of relativity, which is to say, on the preliminary version that preceded his celebrated, ``The field equations of gravitation" of November of that year. Noether must have attended these lectures. 
It is clear from Hilbert's letter
to Einstein of 27 May 1916 
that she had by then already written
some notes on the subject of the problems arising in the general theory of relativity: 
\begin{quotation}
\noindent My law [of conservation] of energy
is probably linked to
yours; I have already given Miss Noether this question to study. 
\end{quotation} 
Hilbert adds that, to avoid a long explanation, he has appended to his letter
``the enclosed note of Miss
Noether''. On 30 May 1916, Einstein answered him
in a brief letter in which he derived a consequence of the equation
that Hilbert had proposed ``which
deprives the theorem of its sense'', and then asks, 
``How can this be
clarified?'' and continues, 
\begin{quotation}
\noindent Of course it would be sufficient if you asked Miss
Noether to clarify this for
me.\footnote{Einstein, {\it Collected Papers}, 8A, nos. 222 and 223.}
\end{quotation}
Thus, her expertise was conceded by both Hilbert and
Einstein as early as her first year in G\"ottingen, and was later acknowledged more explicitly by Klein when he re-published the articles that had appeared in the {\it G\"ottinger Nachrichten} of 1917 and 1918 in his collected works in 1921, a few years before his death.
 

\subsection*{Noether's article of 1918}

In early 1918, Noether published an article on the problem of the invariants
of differential equations in the \selectlanguage{german}{\it G{\"o}ttinger Nachrichten},
``Invarianten beliebiger Differentialausdr\"ucke''\selectlanguage{english} (Invariants of arbitrary differential expressions), 
which was presented by Klein at the meeting of the \selectlanguage{german}{\emph{K\"onigliche Gesellschaft der
    Wissenschaften zu G\"ottingen}}\selectlanguage{english} (Royal G{\"o}ttingen Scientific Society) of 25
January.
It was then, in the winter and spring
of 1918, that Noether discovered the profound reason for the
difficulties that had arisen in the interpretation of the
conservation laws in the general theory of relativity.
Because she had left G\"ottingen for a visit to Erlangen to see her widowed, ailing father, her correspondence remains and it yields an account of her progress in this search. 
In her postcard to Klein of {15 February}, she already sketched her second theorem, but only in a particular case.
It is in her letter to Klein of 
{12 March} that Noether gave a preliminary formulation of an essential consequence of what would be her second theorem, dealing with the invariance of a variational problem under the action of a group which is a subgroup of an infinite-dimensional group.
On 23 July, she presented her results to the \selectlanguage{german}{\emph{Mathematische 
Gesellschaft zu G\"ottingen}}\selectlanguage{english} (G{\"o}ttingen Mathematical Society).
The article which contains her two theorems is ``Invariante Variationsprobleme''
(Invariant variational
  problems). On 26 July, it was
presented by Klein at the meeting of the more important, because it was not restricted to an audience of pure mathematicians, G{\"o}ttingen Scientific Society, and published in the {\it Nachrichten} (Proceedings) of the Society of 1918, on pages 235--247. A footnote on the first page of her article indicates that ``The definitive version of the manuscript was prepared only 
at the end of September."


\subsection*{What variational problems was Noether considering?}
\begin{quotation}
\noindent We consider variational problems which are invariant under a continuous group (in
the sense of Lie). [...] What follows thus depends upon a combination
of the methods of the formal calculus of variations and of Lie's theory of groups.
\end{quotation}
Noether considers a general $n$-dimensional variational problem of 
order $\kappa$ for an ${\mathbb R}^\mu$-valued function, where $n$, $\kappa$ and $\mu$ are arbitrary integers, defined by an integral,
$$
I = \int \cdots \int f \left( x,u, \frac{\partial u}{\partial x},
  \frac{\partial^{2}u}{\partial x^{2}},\cdots,  \frac{\partial^{\kappa}u}{\partial x^{\kappa}}\right)dx,
$$
where $x = (x_1, \ldots, x_n) = (x_\lambda)$ denote the independent variables, and where
$u~=~(u_1, \ldots, u_\mu)= (u_i)$ are the dependent variables.
In footnotes she states her conventions and explains her abbreviated notations:
``I omit the indices here, and in
  the summations as well whenever it is possible, and I write
  $\displaystyle\frac{\partial^{2}u}{\partial x^{2}}$ 
for  $\displaystyle\frac{\partial^{2}u_{\alpha}}{\partial x_{\beta} \partial
    x_{\gamma}}$, etc.'' and ``I write 
$dx$ for $dx_{1}\ldots dx_{n}$ for short.''

Noether then states her two theorems:\footnote{I cite the English translation of Noether's article that appeared in {\it The Noether Theorems} [2010].}
\begin{quotation} 
\noindent In what follows we shall examine the following two
theorems:

\noindent {\bf I.} If the integral $I$ is invariant under a [group] ${\mathfrak{G}}_\rho$, then
there are $\rho$ linearly independent combinations among the
Lagrangian expressions which become divergences -- and
conversely, that implies the invariance of $I$ under a
[group] ${\mathfrak{G}}_\rho$. The theorem remains valid in the limiting
case of an infinite number of parameters.

\noindent {\bf II.} If the integral $I$ is invariant under a [group]
${\mathfrak{G}}_{\infty \rho}$
depending upon arbitrary functions and their derivatives
up to order $\sigma$, then there are $\rho$ identities among the
Lagrangian expressions and their derivatives up to
order~$\sigma$. Here as well the converse is valid.\footnote{In a footnote, Noether announces that she will comment on ``some   trivial exceptions'' in the next section of her article.}
\end{quotation}
Noether  proves the direct part of both theorems in Section 2, then the converse of theorem I in Section 3 and that of theorem II in Section 4.
In Section 2, she assumes that the action integral 
$I = \int f dx$ is invariant. Actually,
she assumes a more restrictive hypothesis, the invariance of the
integrand, $f dx$, which 
is to say $\delta(f dx) = 0$. This hypothesis is
expressed by the relation,
$$
\bar \delta f + {\mathrm{Div}}(f \cdot \Delta x) = 0 .
$$
Here ${\mathrm{Div}}$ is the divergence of vector fields, and $\bar \delta f$ is the variation of $f$ induced by the variation
$$
{\bar{\delta}}u_{i}=\Delta u_{i}- \sum \frac{\partial u_{i}}{\partial
x_{\lambda}} \Delta x_{\lambda}.
$$
Thus, Noether 
introduced the evolutionary
representative, $\bar \delta$, of the  
vector field, $\delta$,  and $\bar \delta f$ is the Lie derivative of $f$
in the direction of the vector field $\bar \delta$.
What she introduced, with the notation $\bar \delta$, 
is a generalized vector field, which is not a vector field in the usual sense, on the trivial 
vector bundle $ {\mathbb R}^n
  \times {\mathbb R}^\mu \to {\mathbb R}^n $.
In fact, if 
$$\delta =  \sum_{\lambda=1}^n X^\lambda(x) \frac{\partial}{\partial{x^\lambda}} 
+ \sum_{i=1}^\mu Y^i(x,u) \frac{\partial}{\partial u^i},
$$
then $\bar \delta$ is the vertical generalized vector field 
$$\bar \delta = \sum_{i=1}^\mu \left(Y^i(x,u) - X^\lambda(x) u^i_\lambda \right)\frac{\partial}{\partial u^i},
$$
where $u^i_\lambda= \frac{\partial u^i}{\partial{x^\lambda}}$. It is said to be ``generalized'' because its components depend on the derivatives of the $u^i(x)$. It is said to be ``vertical'' because it contains no terms in 
$\frac{\partial}{\partial{x^\lambda}}$.\footnote{ The evolutionary
representative of an ordinary vector field has also been called the vertical
  representative. Both terms are modern. Noether does not give $\bar \delta$ a name. An arbitrary vertical generalized vector field is written locally, $Z = \sum_{i=1}^\mu Z^i\left(x,u,\frac{\partial u}{\partial x},
  \frac{\partial^{2}u}{\partial
    x^{2}},\cdots\right)\frac{\partial}{\partial u^i}$.}

By integrating by parts,
Noether obtains the identity 
$$
\sum \psi_{i} \,  \bar{\delta} u_{i} = \bar\delta f + \mathrm{Div}\ A,
$$
where the $\psi_i$'s are the ``Lagrangian
expressions'', {\it i.e.}, the components of the Euler--Lagrange
derivative of $f$, 
and $A$ is linear in $\bar \delta u$ and its
derivatives.
In view of the invariance hypothesis which is expressed by
$\displaystyle{
\bar \delta f + {\mathrm{Div}}(f \cdot \Delta x) = 0}$,
this identity can be written

\centerline{$\displaystyle{\sum \psi_{i} \,  \bar{\delta} u_{i}
    =\mathrm{Div}\ B }$, \, \, where \, \, $B = A - f \cdot \Delta X$.}
\noindent Therefore $B$ is a conserved current 
for the Euler--Lagrange equations of $f$ and the proof of the direct part of Theorem I is complete:
the equations  ${\mathrm{Div}\ B =0}$ 
are the  {conservation laws} that
are satisfied when the  Euler--Lagrange equations ${\psi_i=0}$ are satisfied.

Noether then proves the converse of Theorem I. The existence of $\rho$ ``linearly independent divergence
relations'' implies the invariance under a (Lie) group of symmetries of dimension $\rho$, by passing from the infinitesimal
symmetries to invariance under their flows, provided that the vector fields
$\Delta u$ and $\Delta x$ are ordinary vector fields. 
Thus the existence of $\rho$ linearly independent conservation laws yields the infinitesimal invariance of $f$ under a {Lie algebra of 
infinitesimal symmetries} of dimension $\rho$ but, in the general case, these symmetries are generalized vector fields.
Equivalence relations have to be introduced to make these statements precise.

Theorem II deals with a {symmetry group depending on arbitrary functions}---such as the group
of diffeomorphisms of the space-time manifold and, more generally,
the groups of all gauge theories that would be developped, beginning with the article of Chen Ning Yang
and Robert L. Mills, ``Conservation of isotopic spin and isotopic gauge invariance'', in 1954. Noether showed that to such symmetries there correspond
{identities} satisfied by the variational derivatives, and conversely.
The assumption is that 
``the integral $I$ is invariant under a [group] ${\mathfrak
    G}_{\infty \rho}$ depending upon arbitrary functions and their
 derivatives up to order $\sigma$'', {\it i.e.}, 
Noether 
assumes the existence of $\rho$ infinitesimal symmetries of the variational integral, each of
which depends linearly on an arbitrary function $p$ (depending on
$\lambda
= 1, 2, \ldots, \rho)$ of the variables 
$x_1, x_2, \ldots, x_n$, and its derivatives up to order
$\sigma$. 
Such a symmetry is defined by a {vector-valued linear
differential operator}, ${\mathcal D}$, of order $\sigma$, 
 with components ${\mathcal D}_i$, 
$i = 1, 2, \ldots, \mu$.
Noether then introduces, without giving it a name or a
particular notation, the {adjoint operator},
$({\mathcal D}_i)^*$,
of each of the ${\mathcal D}_i$'s. By construction,
$({\mathcal D}_i)^*$ satisfies
$$\psi_i \, {\mathcal D}_i (p) = 
({\mathcal D}_i)^* (\psi_i)  \, p  \, \, 
+ \, \, 
{\rm Div} \, \Gamma_i,
$$ 
where $\Gamma_i$ 
is linear in $p$ and its derivatives. 
The symmetry assumption and, again, an integration by parts imply
$$
\sum_{i=1}^{\mu}  
\psi_i \, {\mathcal D}_i (p) = {\rm Div} \, B.
$$
This relation implies
$$
\sum_{i=1}^{\mu}  
({\mathcal D}_i)^*(\psi_i) \  p = 
{\rm Div} \large( B - \sum_{i=1}^{\mu}\Gamma_i\large).
$$
Since
$p$ is arbitrary, by Stokes's theorem and the
Du Bois-Reymond lemma, 
$$\sum_{i=1}^{\mu}({\mathcal D}_i)^*(\psi_i)  =0. 
$$
Thus, for each $\lambda = 1, 2, \ldots, \rho$, there is a 
{differential relation}
among the components $\psi_i$ of the Euler--Lagrange
derivative of the Lagrangian $f$ that
is identically satisfied. 

Noether explains 
the precautions that must be taken---the introduction of an
equivalence relation on the symmetries---for 
the converse to be valid. She then observes that each identity may be written
$
\sum_{i=1}^{\mu} a_i \psi_i = {\rm Div} \,
\chi$,
where $\chi$ is defined by a linear 
differential operator acting on the 
$\psi_i$'s. 
She shows that each divergence, ${\rm Div} \, B$, introduced above, is equal to the divergence of a quantity $C$, where $C$ vanishes 
once the Euler--Lagrange
equations, $\psi_i =0$,
are satisfied.
Furthermore, from the equality of the divergences of
$B$ and $C$, it follows that
$$B= C + D
$$
for some $D$ whose divergence vanishes identically, which is
to say, independently of the satisfaction  
of the Euler--Lagrange 
equations.
These are the conservation laws that Noether called {\it improper
 divergence relations}.
In the modern terminology, 
there are two types of trivial conservation laws.
If the quantity $C$ itself, and not only its divergence, vanishes on $\psi_i =0$, then
$C$ is a \emph{trivial conservation law of the first kind}.
If the divergence of $D$ vanishes identically, {\it i.e.}, whether or not $\psi_i =0$,
then $D$ is a  \emph{trivial
  conservation law of the second kind} or ${\rm Div} \,D$ is a \emph{null divergence}.  

\subsection*{Hilbert's conjecture, groups and relativity} 
The last section of Noether's article deals with Hilbert's conjecture. He had asserted, without proof, in early 1918   that, in the case of general relativity, ``the energy equations do not exist at all", that is, there are no proper conservation laws:
\begin{quotation}
\noindent Indeed I claim that for {\it general} relativity, that is, in the case of the {\it general} invariance
of the Hamiltonian function, energy equations which, in your sense, correspond
to the energy equations of the orthogonally invariant theories, do not exist at all; I can even
call this fact a characteristic feature of the general theory of relativity.\footnote{Klein [1917], p. 477, citing Hilbert.}
\end{quotation}
Noether shows that the situation is better understood
``in the
more general setting of group theory''.
She explains the
apparent paradox that arises from the consideration of the
finite-dimensional subgroups of groups 
that depend upon arbitrary functions. She emphasizes the conclusion of her argument by setting it as follows, with italics in the original: 
\begin{quotation}
\noindent {\it Given $I$ 
invariant under the group of translations, then the energy relations
are improper if and only if $I$ is invariant under an
infinite group which contains the group of translations as
a subgroup.} 
\end{quotation}
Noether concludes by quoting in her final footnote Klein's striking formula from page 287 of his 1910 paper, ``\"Uber die geometrischen Grundlagen der Lorentzgruppe'':
\begin{quotation}
\noindent The term
``relativity'' as it is used in physics should be replaced by
``invariance with respect to a group''.
\end{quotation}


\subsection*{How original were Noether's two theorems?}
Noether's article did not appear in a vacuum.
Analysing the contributions of her predecessors requires a detailed development\footnote{See ``{\it The Noether Theorems}", p. 29-39.}. 
Here, I shall only give a very brief account of some of the most important points of this history.
Lagrange, in his {\it M\'echanique analitique}  
(1788), 
claimed that his method for deriving 
``a general formula for the motion of bodies'' yields ``the general equations that contain the principles, or
theorems known under the names of the {\it conservation of kinetic energy}, of the {\it conservation of the motion of the
center of mass}, of the {\it conservation of the momentum of
rotational motion}, of the {\it principle of areas}, and of the
{\it principle of least action}''.\footnote{Lagrange [1788], p. 182, italics in the original.}
In the second edition of his {\it M\'ecanique
analytique}, in 1811, as a preliminary to his treament of dynamics, he presented a detailed history of the diverse ``principes
ou th\'eor\`emes'' (principles
or theorems) formulated before his {\it M\'ecanique}, thus recognizing the contributions of his predecessors in the  discovery of these principles -- 
Galileo, Huyghens, Newton, Daniel Bernoulli, Maupertuis, Euler, the Chevalier Patrick d'Arcy 
 and d'Alembert--, and in this second edition, he explicitly  observed a correlation between these principles of conservation and invariance properties.  
After Lagrange, the correlation between invariances and conserved quantities was surveyed 
by {Jacobi} in several chapters of his 
{{\it Vorlesungen \"uber Dynamik}}, lectures delivered in 1842-43 but only published posthumously in 1866.
The great advances of Sophus Lie (1842--1899), his theory of continuous groups of transformations that was published in articles and books that appeared between 1874 and 1896, became the basis of all later developments, such as the work of Georg Hamel (1877--1954) on the calculus of variations and mechanics in {1904}, and the publication of
{Gustav Herglotz} (1881--1953) on the 10-parameter invariance group of the [special] theory of relativity
in {1911}. 
In her 1918 article, Noether cited Lie very prominently as his name appears three times in the eight lines of the introductory paragraph, but with no precise reference to his published work. Both Hamel and Herglotz were cited by her. 
In her introduction, she also referred to publications, all of them still very recent, by ``[Hendrik] Lorentz and his students (for example, [Adriaan Daniel] Fokker), [Hermann] {Weyl}, and {Klein} for certain infinite groups'' and, in a footnote, she wrote, ``In a paper by [Adolf] Kneser that has just appeared (Math. Zeitschrift, vol. 2), the determination of invariants is dealt with by a similar method.'' In fact, while Noether was completing the definitive version
of her manuscript, in August 1918, Kneser had  
submitted an article, ``Least action and Galilean
relativity", in which he used Lie's infinitesimal
transformations and, as Noether would do, emphasized the relevance
of Klein's Erlangen program, but he did not treat questions of invariance.
Noether stressed the relation of her work to ``Klein's second note, {\it G\"ottinger Nachrichten}, 19 July 1918",\footnote{Klein [1918].} stating that her work and Klein's were ``mutually influential'' and referring to it for a more complete bibliography.
In section 5 of her paper, she cited an article ``On the ten general invariants of classical mechanics'' by Friedrich Engel (1861--1941) that had appeared two years earlier.
Indeed, scattered results in classical and relativistic mechanics, tying
together properties of invariance and conserved quantities, had already appeared
in the publications of Noether's predecessors which she acknowledged.
However, none of  them had discovered
the general principle contained in her Theorem I and its converse.
Her Theorem II and its converse were completely new. 
In the expert opinion of the theoretical physicist Thibaut Damour,\footnote{Damour is a professor at the Institut des Hautes \'Etudes Scientifiques and a member of the Acad\'emie des Sciences de l'Institut de France.} the second theorem should be considered the most important part of her article. It is certainly the most original.

\subsection*{How modern were Noether's two theorems?}
What Noether simply called ``infinitesimal transformations'' are, in fact, vast generalizations of the ordinary vector fields, and are now called ``generalized vector fields''. 
They would eventually be re-discovered, independently, in 1964 by Harold Johnson, who called them ``a new type of vector fields'', and in 1965 by Robert
Hermann. They appeared again in 1972 when Robert L. Anderson, Sukeyuki Kumei and Carl
Wulfman published their ``Generalization of the 
concept of invariance of differential equations. Results of applications to some Schr\"odinger equations" in 
{\it Physical Review Letters}. In 1979, R. L. Anderson, working at the University of Georgia, in the United States, and Nail Ibragimov (1938--2018), then a member of the Institute of Hydronymics in the Siberian branch of the USSR Academy of Sciences in Novosibirsk -- such east-west collaboration was rare at the time --, in their monograph,
{\it Lie-B\"acklund
Transformations in Applications}, duly citing Klein and Noether while claiming to generalize ``Noether's classical theorem'', called them
``Lie-B\"acklund transformations'', a misleading term because Albert~V. B\"acklund (1845--1922) did not introduce this vast generalization of the concept of vector fields, only infinitesimal contact transformations. 
The concept of a generalized vector field is essential in the theory of integrable systems which became the subject of intense research after 1970. On this topic,
Noether's work was modern, half-a-century in advance of these re-discoveries.
Peter Olver's book [1986] is both a comprehensive handbook of the theory of generalized symmetries of differential and partial differential equations, and the reference for their history, while his article of the same year on ``Noether's theorems and systems of Cauchy--Kovalevskaya type'' is an in-depth study of the mathematics of Noether's second theorem. 
His article [2018], written for the centenary of Noether's article, stresses the importance of her invention of the generalized vector fields.

In G\"ottingen, Noether had only one immediate follower,
{Erich Bessel-Hagen} 
(1898--1946), who was Klein's student.
In 1921, he published an article in the {\it Mathematische Annalen},  entitled \selectlanguage{german}{\it \"Uber die Erhaltungss\"atze der Elektrodynamik}\selectlanguage{english} (On the conservation laws of
electrodynamics), in which he
determined in particular those conservation laws that are the result of the
conformal invariance of Maxwell's equations. In this paper, Bessel-Hagen recalls that it was Klein who had posed
the problem of ``the application to Maxwell's equations 
of the
theorems stated by Miss Emmy Noether} about
two years ago
regarding the invariant variational problems"
and he writes that, in the present paper, he formulates the two Noether theorems ``slightly more
generally'' than they had been formulated in her article.
How did he achieve this more general result? By introducing the concept of ``divergence symmetries'' which are infinitesimal transformations which leave the Lagrangian invariant up to a divergence term, or ``symmetries up to divergence''.
They correspond, not to the
invariance of the Lagrangian $f dx$, 
but to the invariance of the action integral
$\int f dx$,
 {\it i.e.}, 
instead of satisfying the condition
$\delta (f dx)= 0$, they satisfy, 
the weaker condition
$\delta (f dx)=  \mathrm{Div} C$,
where $C$ is a vectorial expression.
Noether's fundamental relation remains valid under this
weaker assumption, provided that  
$B = A - f\cdot \Delta x$
be replaced by 
$B = A + C - f\cdot \Delta x.$
Immediately after he stated that he had proved the theorems in a slightly more general form than Noether had, Bessel-Hagen added: ``{I owe these [generalized theorems] to an oral communication by Miss Emmy Noether herself}''. We infer that, in fact,
this more general type of symmetry was also Noether's invention. Bessel-Hagen's acknowledgment is evidence that, to the question, ``Who invented divergence symmetries?'', the answer is ``Noether''. 

\subsection*{How influential were Noether's two theorems?}

The history of the {reception} of Noether's article in the years 1918--1970 is surprising.
She submitted the ``Invariante Variationsprobleme'' for her
{\it Habilitation}, finally obtained in 1919,
but she never referred to her article in any of  her subsequent publications. I know of only one mention of her work of 1918 in her later writings, in a letter she sent eight years later to Einstein who was then an editor of the journal {\it Mathematische Annalen}. In this letter, which is an informal referee report, she rejects a submission ``which unfortunately is by no means suitable" for the journal, on the grounds that ``it is first of all a restatement that is not at all clear of the principal theorems of
my \selectlanguage{german}`Invariante Variationsprobleme'\selectlanguage{english} (G\"ott[inger] Nachr[ichten], 1918 or 19), with a slight generalization--the invariance of the integral up to a divergence term--which
can actually already be found in Bessel-Hagen (Math[ematische] Ann[alen], around
1922).''\footnote{For a facsimile, a transcription, and a translation of Noether's letter, see {\it The Noether Theorems} [2010], p. 161--165, and see comments on this letter, {\it ibid.}, p. 51--52.}

I found very few early occurrences of Noether's title in books and articles. 
While {Hermann Weyl}, in {\it Raum, Zeit, Materie}, first published in 1918,
performed computations
very similar to hers, he referred to Noether only once, in a
footnote in the third (1919) and subsequent editions. It is clear that {Richard Courant} must have been aware of her work because a brief
summary of a limited form of both theorems appears in all
German, and later English editions of ``Courant--Hilbert'', the widely read treatise on methods of mathematical physics
first published in 1924.
It is remarkable that we found so few explicit mentions of Noether's results in searching the literature of the 1930's.
In 1936, the little known physicist, {Moisei A. Markow} 
(1908--1994), who was a member of the Physics Institute of the
U.S.S.R. Academy of Sciences in Moscow, published
an article in the {\it Physikalische Zeitschrift der
Sowjetunion} in which he refers to ``the well-known theorems of
Noether.''
Markow was a former student of 
Georg 
B. Rumer (1901--1985) who had been an assistant of Max
Born in G\"ottingen from 1929 to 1932. Rumer, in 1931, had
proved the Lorentz invariance of the Dirac operator
but did not allude to any associated conservation laws, while in his
articles on the general theory of relativity published in the {\it G\"ottinger
Nachrichten} in 1929 and 1931, he cited Weyl but never Noether.
Similarly, it seems that V.~A. Fock (1898--1974) never referred to Noether's work in any of his papers to which it was clearly relevant, such as his celebrated ``Zur Theorie des Wasserstoffatoms'' (On the theory of the hydrogen atom) of 1935. 
Was it because, at the time, papers carried few or no citations? or because Noether's results were considered to be ``classical''? The answers to both questions are probably positive, this paucity of citations being due to several factors.         

An early, explicit reference to Noether's publication is found in the article of Ryoyu Utiyama (1916--1990), then in the department of physics of Osaka Imperial University, ``On the interaction of mesons with the gravitational field.~I.'', which appeared in {\it Progress of Theoretical Physics} [1947], four years before he was awarded the Ph.D. His paragraph I begins with the ``Theory of invariant variation'' for which he cites both Noether's 1918 article and p. 617 of Pauli's ``Relativit\"atstheorie'' [1921]. Following Noether closely, he proves the first theorem, introducing ``the substantial variation of any field quantity'', which he denotes by $\delta^*$, {\it i.e.}, what Noether had denoted by $\bar \delta$, and also treats the case where the dependent variables ``are not completely determined by [the] field equations but contain $r$ undetermined functions''. This text dates, in fact, to 1941 as the author reveals in a footnote on the first page: ``This paper was published at the meeting[s] of [the] Physico-mathematical Society of Japan in April 1941 and October 1942, but because of the war the printing was delayed''. Such a long delay in the publication of this scientific paper is one example -- among many -- of the influence of world affairs on science. It appears that this publication is a link in the chain leading from Noether's theorems to the development, by the physicists, of the gauge theories, where the variations of the field variables depend on arbitrary functions. Episodes in this history, told by Utiyama himself, were published in Lochlainn O'Raifeartaigh's book [1997], from which we learn that, although Utiyama published his important paper ``Invariant theoretical interpretation of interaction'' only in 1956, two years after the famous article of Yang and Mills, he had worked independently and had treated more general cases, showing that gauge potentials are in fact affine connections. In this paper, Utiyama gave only six references: one is (necessarily) to the publication of Yang and Mills, another is to his own 1947 paper, clearly establisihing the link from his previous work to the present one, and another reference is to p. 621 of Pauli [1921]. This time, however, a reference to Pauli serves as a reference to Noether, so that her name does not appear.

In later developments, in the Soviet Union
in 1959, Lev S. Polak published a translation of Noether's 1918
article into {Russian} and, in 1972, Vladimir Vizgin published a historical monograph whose title, in English translation, is {\it  The Development of the interconnection between invariance principles and conservation laws in classical physics}, in which he analyzed both 
Noether's theorems.
At that time, new formulations of Noether's first theorem had started to appear with the textbook of Israel M. Gelfand and Sergej V. Fomin on the
calculus of variations, published in Moscow in 1961, which contains a modern presentation of Noether's first
theorem -- although not yet using the formalism of jets as would soon be the case --, followed by a few lines about her second theorem. This text appeared in an English translation two years later.
In the 1970s, Gelfand published several articles with Mikhael Shubin, {Leonid Dikii} (Dickey), {Irene Dorfman},
and Yuri Manin on the ``formal calculus of variations'', not mentioning Noether because they dealt mainly with the Hamiltonian formulation of the problems, while Manin's ``Algebraic
theory of nonlinear differential equations'' [1978] as well as  
{Boris Kupershmidt}'s ``Geometry of jet bundles and the structure of Lagrangian and Hamiltonian formalisms" [1980] both contain a ``formal Noether theorem'', which is a modern, generalized version of her first theorem. 
A few years earlier, in the article, ``Lagrangian formalism in the calculus of variations'' [1976], Kupershmidt had already presented an invariant  approach  to  the  calculus  of variations  in  differentiable  fibre  bundles and formulated     Noether's first  theorem  for the Lagrangians of arbitrary finite order.

Further research in geometry in Russia yielded new genuine 
generalizations of the concepts introduced by Noether and of her results.  
Alexandre Vinogradov (1938--2019), who had been a member of
Gelfand's seminar in Moscow, left the Soviet Union for Italy in 1990 and the second part 
of his career was at the University of Salerno.
Beginning in 1975, together with {Joseph   
Krasil'shchik}, who worked in Moscow, for several years in the Netherlands, and again in Moscow at the Independent University, Vinogradov
published extensively on symmetries, at a very general and
abstract level, greatly generalizing Noether's formalism and results, and on their
applications, a theory fully expounded in their book of 1997. 

Searching for other lines of transmission of Noether's results, one finds that, in the early 1960s, {Enzo Tonti} (later professor at the University of
Trieste) translated Noether's article into Italian but his translation has
remained in manuscript. It was transmitted to {Franco Magri} 
in Milan who, in 1978, wrote an article in Italian
where he clearly set out the relation between symmetries and
conservation laws for non-variational equations, a significant development, but he did not treat the case of operators defined on manifolds.

In France, Jean-Marie Souriau (1922--2012), was well aware
of ``les m\'ethodes d'Emmy Noether'' which he cited  as early as 1964, on page 328 of his first book, \selectlanguage{french}{\it G\'eom\'etrie et relativit\'e}\selectlanguage{english}.
In 1970, independently of
Bertram Kostant (1928--2017), he introduced the concept of a moment map.
The conservation of the moment of a Hamiltonian action is the Hamiltonian version of Noether's first theorem. Souriau called that result 
{``le th\'eor\`eme de Noether symplectique''} although there is nothing Hamiltonian or symplectic in Noether's article!
Souriau's fundamental work on symplectic geometry and mechanics was based on Lagrange, as he himself claimed, and it was also a continuation of Noether's theory. 

\subsection*{From general relativity to cohomological physics}
The history of the second theorem -- the improper conservation laws --
is part of the history of general relativity. In the literature on                                                            the general theory, 
the improper conservation laws which are ``trivial of
the second kind'' are called ``strong laws'', while the conservation laws
obtained from the first theorem are called 
``weak laws''. 
The strong laws play an important role 
in basic papers of 
Peter~G. Bergmann in 1958, of
Andrzej Trautman in 1962, and
of Joshua~N. Goldberg in 1980. While the second theorem, which explained in which cases such improper conservation laws would exist, had been known among relativists since the early 1950s, it became an essential tool in the non-abelian gauge theories that were developed by the mathematical physicists, following the publication of Utiyama's paper in 1956 that generalized the Abelian theory of Yang and Mills of 1954. 

The identities that were proved by Noether in her second theorem are at the
basis of Jim Stasheff's ``cohomological physics''. They appeared already in his lecture at Ascona [1997]. Then, in his article with Tom Fulp and Ron Lada published in 2003,  
``Noether's variational theorem II and the BV formalism'', Noether's identities associated with the infinitesimal gauge symmetries
of a Lagrangian theory appear as the anti-ghosts in the Batalin--Vilkovisky construction for the quantization
 of Lagrangians with symmetries.
The validity of Noether's second theorem is extended to ever more general kinds of symmetries, interpreting physicists' constructions in gauge theories of increasing complexity.

\subsection*{Have the Noether theorems been generalized?}
Whether the Noether theorems have been generalized has a straightforward answer: except for Bessel-Hagen (and we have seen that his generalization was certainly suggested and  probably entirely worked out by Noether herself), it was not until the 1970s.
Until then, the {so-called ``generalizations''}
were all due to physicists and mathematicians who had no direct knowledge of her article but still thought that they were generalizing it, while they were generalizing 
the truncated and restricted version of her first theorem that had appeared in
{Edward L. Hill}'s article,
``Hamilton's principle and the conservation theorems of
mathematical physics'', in 1951.

In the late 1970's and early 1980's, using different languages, both linguistically and mathematically, Olver, in Minneapolis, and Vinogradov, in Moscow, made great advances in the Noether theory. 
Equivalence classes were defined for symmetries on the one hand and for conservation laws on the other, bringing precision to the formulation of Noether's results.
In order to set up a one-to-one map between symmetries and conservation laws it is appropriate to first consider 
the enlarged class of the divergence symmetries (which are the infinitesimal tranformations leaving the 
Lagrangian invariant up to a divergence term, {\it i.e.}, the concept of symmetry to be found in Bessel-Hagen's article of 1921). Define a divergence symmetry of a differential equation to be trivial 
if its evolutionary representative 
vanishes on the solutions of the equation or if its divergence vanishes identically, independently of the field equations,
and consider the equivalence classes of divergence symmetries modulo the trivial symmetries.
Recall the definition of the trivial conservations laws of the first and of the second kind and consider the equivalence classes of conservation laws modulo the trivial ones.
Restrict the consideration of Lagrangians to those whose system of Euler--Lagrange equations is ``normal'', meaning roughly that the highest-order partial derivatives of the unknown functions are expressed in terms of all the other derivatives.
Then one can formulate what can be called ``the Noether-Olver-Vinogradov theorem" which took the following form, both rigorous and simple, ca. 1985:
\begin{quotation}
\noindent For Lagrangians such that the Euler--Lagrange equations are a normal
system, Noether's correspondence induces a one-to-one map 
between equivalence classes of divergence symmetries 
and equivalence classes of conservation laws.
\end{quotation}

Concerning the extension to non-variational equations of Noether's correspondence between symmetries and conservation laws, we find the early work of 
Magri [1978] who showed that, 
if ${\mathcal D}$ is a differential operator and $V{\mathcal D}$ is its linearization,
searching 
for the restriction of the
kernel of the adjoint $(V{\mathcal D})^*$ of $V{\mathcal D}$
to the solutions of ${\mathcal D}(u) =0$ is an algorithmic method for
the determination of the
conservation laws for a possibly non-variational equation, ${\mathcal D}(u) = 0$.
For an Euler--Lagrange operator, the linearized operator is  self-adjoint. Therefore 
{this result generalizes the first Noether theorem.}
This idea is to be found later and much developped in the work of several mathematicians, 
most notably Vinogradov, 
Toru Tsujishita, Ian Anderson, and Olver.


Meanwhile, Noether's theory was being set in the modern language of differential
geometry and generalized.
Trautman, in his ``Noether equations and conservation laws'' [1967], followed by ``Invariance in Lagrangian systems'' [1972], was the first to present even a part of Noether's article in the language of manifolds, fiber 
bundles and,
in particular, the jet bundles that had been defined and studied around 1940 by
Charles Ehresmann (1905--1979) and his student Jacques Feldbau (1914--1945), and by Norman Steenrod (1910--1971).
In 1970, Stephen Smale published the first part of his article on ``Topology and mechanics'' in which he proposed a {geometric framework for mechanics on the tangent bundle} of a manifold.
Hubert Goldschmidt and
Shlomo Sternberg 
wrote a landmark paper in 1973
in which they formulated the Noether theory for first-order Lagrangians in an {intrinsic, geometric fashion.}
Jerrold Marsden published extensively on the {theory and important applications} of Noether's correspondence between invariance and conservation from 1974 until his death at the early age of 68 in 2010. 
In the 1970s, several other authors contributed to the ``geometrization'' of Noether's first theorem, notably Jedrzej \'Sniatycki, Demeter Krupka, and Pedro Garc{\'i}a.
The ideas that permitted the recasting of Noether's theorems in
geometric form and their genuine generalization were first of all
that of {smooth differentiable manifold} ({\it i.e.}, manifolds of class $C^\infty$), and then
the concept of a {jet of order $k$} of a mapping ($k \geq 0$) defined as the collection of the values of the components
in a local system of coordinates of a vector-valued function and of
their partial derivatives up to the order $k$, the concept of
manifolds of jets of sections of a fiber bundle, and finally of jets of
infinite order. The {manifold of jets of infinite order of sections} of a
fiber bundle is not defined directly but as the
{inverse limit} 
of the manifold of jets of order $k$, when
$k$ tends to infinity. 
It was {Vinogradov} who showed in 1977 that the generalized vector
fields are nothing other than ordinary vector fields on the bundle
of jets of infinite order of sections of a bundle. 
Both Lagrangians and conservation laws then appear as special types
of differential forms. 
The divergence operator may be interpreted as a
horizontal differential, one that acts on the independent variables
only, and that yields a {cohomological interpretation} of
Noether's first theorem.
The study of the {exact sequence of the
  calculus of variations}, and of the
{variational bicomplex}, which constitutes a vast
generalization of Noether's theory,
 was
developed in 1975 and later by {W{\l}odzimierz Tulczyjew} in Warsaw, {Paul Dedecker} in Belgium, {Vinogradov} in Moscow, Tsujishita in Japan, and, in the United States, by Olver
and by {Ian Anderson}.



In the discrete versions of the Noether theorems, the differentiation operation is replaced by a shift operator.
The independent variables are now integers, and the integral is replaced by a sum, ${\mathcal L}[u]= \sum_n L(n,[u])$,
where $[u]$ denotes $u(n)$ and finitely many of its shifts.
The variational derivative is expressed in terms of the inverse shift.
A pioneer was
John David Logan who published ``First integrals in the discrete variational
calculus'' in 1973. Much more recent advances on the {discrete analogues} of the Noether theorems, an active and important field of research, may be found in a series of papers by {Peter Hydon} and {Elizabeth Mansfield}, published since 2001, including a discrete version of the second theorem [2011].

\subsection*{Were the Noether theorems ever famous?}
Whereas both of Noether's theorems were analyzed by Vizgin in his 1972 monograph on invariance principles and conservation laws in classical physics, it appears that, except for an article written that year by Logan where restricted versions of each of her theorems are indeed formulated, the existence of the first and second theorem in one and the same publication was not expressed in written form in any language other than Russian before 1986, when the first edition of Olver's book and his  article where ``Noether's theorems" appear in the title were published.
At roughly the same time, one can find ``theorems'', in the plural, in 
a few other publications: in Hans A. Kastrup's contribution to {\it Symmetries in Physics (1600--1980)}, the proceedings of a meeting held in 1983 in Sant Feliu de Gu\'ixols in Catalonia, published in 1987,
and in my paper, ``Sur les th\'eor\`emes de Noether'', in the proceedings of the ``Journ\'ees relativistes" organized by 
Andr\'e Lichnerowicz in Marseille--Luminy in 1985 which also appeared in 1987.

Fame came eventually.
I quote from Gregg Zuckerman's
``Action principles and global geometry'' [1987]:
\begin{quotation}
\noindent E. Noether's famous 1918
paper, ``Invariant variational problems'' crystallized 
essential mathematical relationships among symmetries, conservation
laws, and identities for the variational or `action' principles
of physics. [\dots]  Thus, Noether's abstract analysis continues
to be relevant to contemporary physics, as well as to applied
mathematics.\footnote{Here Zuckerman cites Olver's {\it Applications of Lie Groups to Differential Equations}.}
\end{quotation}
\noindent Therefore, approximately seventy years after her article had appeared in
the {\it G\"ottinger Nachrichten}, fame came to Noether for this (very small) part of her mathematical {\it {\oe}uvre}.
In 1999, in  the twenty-page contribution of Pierre Deligne and Daniel Freed to the monumental treatise,
{\it Quantum Fields and Strings: A Course for Mathematicians}, she was credited, not only with ``the Noether theorems'', but also with ``Noether charges'' and ``Noether currents''. For as long as gauge theories had been developing, these terms had, in fact, been in the vocabulary of the physicists, such as Utiyama, Yuval Ne'eman, or Stanley Deser whose discussion of ``the conflicting role of Noether's two great theorems" and ``the physical impact of Noether's theorems'' continues to this day in papers and preprints.
At the end of the twentieth century, the importance of the concepts she had introduced 
was finally recognized and her name was attached to them by mathematicians and mathematical physicists alike.

\subsection*{In lieu of conclusion}

One can read in a text published as late as 2003 by a well-known philosopher of science that ``Noether's theorems can be generalized to handle transformations that depend on the $u^{(n)}$ as well.'' 
Any author who had only glanced at Noether's paper, or read parts of Olver's book, would have been aware that Noether had already proved her theorems under that generalized assumption.
This, in fact, is one of the striking and important features of Noether's 1918 article.
Therefore, {\it caveat lector}! It is better to read the original than to rely on second-hand accounts.
For my part, I shall not attempt to draw any philosophical conclusions from what I have sketched here of Noether's ``Invariante Variationsprobleme'', its genesis, its consquences and its influence, because I want to avoid the mistakes of a non-philosopher, of the kind that amateurs make in all fields.

It is clear that Noether was not a proto-feminist. She was not a practicing Jew. Together with her father, she converted to protestantism in 1920, which did not protect her from eventual dismissal from the University of G\"ottingen by the Nazis. She was not an admirer of the American democracy and her sympathies were with the Soviet Union. 
Even though her work of the year 1918 was clearly inspired by a problem in physics, 
she was never herself a physicist and did not return to physics in any of her subsequent publications.
She never explored the philosophical underpinnings or outcomes of her work, in a word, she was not a philosopher. 
She was a generous woman admired by her colleagues and students, and
a great mathematician.

While the Noether theorems derive from the algebraic theory of invariants developed in the nineteenth-century -- a chapter in the history of pure mathematics --, it is clear from the testimony of Noether herself that the immediate motivation for her research was a question that arose in physics, at the time when the general theory of relativity was emerging -- a fact that she stated explicitly in her 1918 article. The results of this article have indeed become -- in increasingly diverse ways which deserve to be much more fully investigated than time and space permitted -- a fundamental instrument for mathematical physicists. On the one hand, these results are essential parts of the theories of mechanics and field theory and many other domains of physical science, and on the other, in a series of mainly separate developments, her results have been generalized by pure mathematicians to highly abstract levels, but that was not accomplished in her lifetime. Had she lived longer, 
she would have witnessed this evolution and the separate, then re-unified, paths of mathematics and physics,
and we are free to imagine that she would have taken part in the mathematical discoveries that issued from her 
twenty-three-page article.

\medskip

\medskip

\subsection*{References}

\m 

\m

\noindent{\sc Anderson} (Robert L.) and {\sc Ibragimov} (Nail H.)

[1979] {\it Lie-B\"acklund
Transformations in Applications}, SIAM Studies in Applied Mathematics,
Philadelphia: SIAM, 1979.

\medskip

\noindent {\sc Bessel-Hagen} (Erich)

[1921] ``\"Uber die Erhaltungss\"atze der Elektrodynamik'', {\it Mathematische Annalen}, 84 (1921),
p. 258--276.
\m

\noindent{\sc Brading} (Katherine A.) and {\sc Brown} (Harvey R.)

[2003] ``Symmetries and Noether's theorems'', in 
{\it Symmetries in Physics, Philosophical Reflections}, Katherine Brading and Elena Castellani, eds., Cambridge: Cambridge University Press, 2003, p. 89--109.

\medskip

\noindent{\sc Brewer} (James W.) and {\sc Smith} (Martha K.), eds.

[1981] {\it Emmy Noether, A Tribute to Her Life and Work}, Monographs and Textbooks in Pure and
Applied Mathematics, vol. 69, New York, Basel: Marcel Dekker, 1981.

\medskip

\noindent{\sc Deligne} (Pierre) and  {\sc Freed} (Daniel S.)

[1999] ``Classical field theory'', in {\it Quantum Fields and Strings: A Course for Mathematicians},
Pierre Deligne, Pavel Etingof, Daniel S. Freed, Lisa C. Jeffrey, David Kazhdan, John W. Morgan,
David R. Morrison, Edward Witten, eds., vol. 1, Providence, RI: American Mathematical
Society/Institute for Advanced Study, 1999, p. 137--225.

\medskip

\noindent{\sc Dick} (Auguste)

[1970] {\it Emmy Noether 1882--1935}, Beihefte zur Zeitschrift
``Elemente der Mathematik", no. 13,
Basel: Birkh\"auser, 1970; English translation by Heidi I. Blocher, Boston, Basel,
Stuttgart: Birkh\"auser, 1981.

\m

\noindent{\sc Euler} (Leonhard) 

[1744]
{\it Methodus inveniendi lineas curvas, maximi minimive proprietate
  gaudentes}, Lausanne and Geneva, {1744}, p. 245.

\medskip

\noindent{\sc Fulp} (Ronald O.), {\sc Lada} (Thomas J.), {\sc Stasheff} (James D.)

[2003] ``Noether's variational theorem II and the BV formalism'' (Proceedings of the 22nd Winter
School ``Geometry and Physics'', Srn{\'i}, 2002), {\it Rendiconti del Circolo Matematico di Palermo},
Series II, Supplement no. 71 (2003), p.~115--126.

\m

\noindent{\sc Gelfand} (Israel M.) and {\sc Fomin} (Sergej V.)

[1961] {\it Variacionnoe ischislenie} [Russian], Moscow: Gosudarstvennoe Izdatelstsvo fizikomatematicheskoj
literatury, 1961; English translation by Richard A. Silverman, Calculus of Variations, Englewood Cliffs, NJ: Prentice-Hall, 1963.

\m

\noindent{\sc Hydon} (Peter) and {\sc Mansfield} (Elizabeth)

[2011] ``Extensions of Noether's second theorem: from continuous to discrete systems'',
{\it  Proc.  R. Soc. Lond. Ser. A Math. Phys. Eng. Sci.} 467 (2011), no. 2135, p. 3206--3221.

\medskip

\noindent{\sc Kastrup} (Hans A.)

[1987] ``The contributions of 
Emmy Noether, Felix Klein and Sophus Lie to the modern concept 
of symmetries in physical systems'', in {\it Symmetries in Physics (1600--1980)},
1st International Meeting on the History of Scientific
Ideas, Sant Feliu de Gu{\'i}xols, Catalonia, Spain (1983), 
 Manuel G. Doncel,
Armin  Hermann, Louis Michel, Abraham Pais, eds.,
Bellaterra (Barcelona): Seminari
d'Hist{\`o}ria de les Ci{\`e}ncies, Universitat Aut{\`o}noma de Barcelona, 1987, p. 113--163.

\m

\noindent {\sc Kimberling} (Clark H.)

[1972] ``Emmy Noether'', {\it American Mathematical Monthly}, 79, no. 2 (1972), p. 136--149. Addendum,
{\it ibid.}, 79, no. 7 (1972), p. 755.

[1981] ``Emmy Noether and her influence'', in Brewer and Smith, eds., 1981, p. 3--61.

\m

\noindent{\sc Klein} (Felix)

[1910] ``\"Uber die geometrischen Grundlagen der Lorentzgruppe'', {\it Jahresbericht der Deutschen
Mathematiker-Vereinigung}, 19 (1910), p. 281--300; reprint, {\it Physikalische Zeitschrift}, 12 (1911),
p. 17--27; {\it Gesammelte mathematische Abhandlungen}, vol. 1, p. 533--552.

[1917] ``Zu Hilberts erster Note \"uber die Grundlagen der Physik'', 
{\it Nachrichten von der K\"oniglichen Gesellschaft der Wissenschaften zu G\"ottingen, Mathematisch-physikalische Klasse} (1917), p. 469--482 (meeting of 25 January 1918); {\it Gesammelte mathematische Abhandlungen}, vol. 1,
p. 553--567.

[1918] ``\"Uber die Differentialgesetze f\"ur die Erhaltung von Impuls und Energie in der Einsteinschen
Gravitationstheorie", {\it Nachrichten von der K\"oniglichen Gesellschaft der Wissenschaften zu G\"ottingen, Mathematisch-physikalische Klasse} (1918), p. 171--189 (meeting of 19 July 1918);
{\it Gesammelte mathematische Abhandlungen}, vol. 1, p. 568--585.

[1921] {\it Gesammelte mathematische Abhandlungen}, 3 vol., Berlin: J.~Springer, 1921--1923; reprint,
Berlin, Heidelberg, New York: Springer, 1973.

\m

\noindent{\sc Kosmann-Schwarzbach} (Yvette)

[1987] 
``Sur les th\'eor\`emes de Noether'', in {\it G\'eom\'etrie et physique} (Journ\'ees relativistes de
Marseille-Luminy, 1985), Y. Choquet-Bruhat, B. Coll, R. Kerner, A. Lichnerowicz, eds.,
Travaux en Cours, 21, Paris: Hermann, 1987, p.~147--160.

[2004]
{\it Les th\'eor\`emes de Noether, Invariance et lois de conservation au XX$^e$ si\`ecle}, \'Editions de l'\'Ecole polytechnique, 2004; revised edition 2006.

[2010]
{\it The Noether Theorems, Invariance and Conservation Laws in the Twentieth Century}, translated by Bertram E.~Schwarzbach,  
Springer, 2010; corrected publication 2018.

\m

\noindent{\sc Krasil'shchik} (Joseph S.) and {\sc Vinogradov} (Alexandre M.), eds.

[1997] {\it Simmetrii i zakony sokhranenija uravnenij matematicheskoj fiziki} [Russian], Moscow:
Factorial, 1997; English translation by A. M. Verbovetsky and I.~S. Krasil'shchik, {\it Symmetries and
Conservation Laws for Differential Equations of Mathematical Physics}, Providence, RI: American
Mathematical Society, 1999.

\m

\noindent{\sc Kupershmidt} (Boris A.)

[1976] ``Lagrangian formalism in the calculus of variations''  [Russian], {\it  Funk\-cional. Anal. i Prilo{\v z}en.} 10 (1976), no. 2, p. 77--78.

[1980] ``Geometry of jet bundles and the structure of Lagrangian and Hamiltonian formalisms'', in
Geometric Methods in Mathematical Physics (Lowell, MA, 1979), Lecture Notes in Mathematics,
775, G. Kaiser and J. E. Marsden, eds., Berlin, Heidelberg, New York: Springer-Verlag, 1980.

\m

\noindent{\sc Lagrange} (Joseph Louis)

[1788] {\it M\'echanique analitique}, Paris: chez la Veuve Desaint, 1788.

[1811] {\it M\'ecanique analytique}, 2nd ed., Paris: Mme Veuve Courcier, vol. 1, 1811.

\m

\noindent{\sc Logan} (John David)

[1973] ``First integrals in the discrete variational
calculus'', {\it  Aequationes
Mathematicae}, 9 (1973), p.~210--220.


\newpage

\noindent{\sc Magri} (Franco)

[1978] ``Sul legame tra simmetrie e leggi di conservazione nelle teorie di campo classiche'', 4th
National Congress of Theoretical and Applied Mechanics, AIMETA, Florence, 25--28 October
1978.

\m

\noindent{\sc Manin} (Yuri I.)

[1978] ``Algebraic
theory of nonlinear differential equations''  
[in Russian], {\it Itogi Nauki i Tekhniki,
Sovremennye Problemy Matematiki}, 11 (1978), p.~5--152; English translation, {\it Journal of Soviet
Mathematics}, 11 (1979), p. 1--122.

\medskip

\noindent{\sc Noether} (Emmy)

[1907] ``\"Uber die Bildung des Formensystems der tern\"aren biquadratischen Form'', {\it Sitzungsberichte der Physikalisch-medizinischen Societ\"at zu Erlangen}, 39 (1907), p. 176--179; {\it Gesammelte Abhandlungen}, p. 27--30.

[1908] ``\"Uber die Bildung des Formensystems der tern\"aren biquadratischen Form'', {\it Journal f\"ur
die reine und angewandte Mathematik}, 134 (1908), p. 23--90, with two tables; {\it Gesammelte  Abhandlungen},
p. 31--99.

[1918] 
\selectlanguage{german}``Invariante Variationsprobleme'', {\it Nachrichten
von der K\"oniglichen Gesellschaft der Wissenschaften zu G\"ottingen,
Mathematisch-physikalische Klasse}\selectlanguage{english} (1918), 
p.~235--257; 
{\it Gesammelte Abhandlungen}, p. 248--270. English translation,
``Invariant variation problems'', by M. A. Tavel, {\it Transport
Theory and Statistical Physics}, 1 (3) (1971), p. 186--207. English translation,
``Invariant variational problems", in Kosmann-Schwarzbach [2010], p. 1--22.

[1983] {\it Gesammelte Abhandlungen/Collected papers},
N. Jacobson, ed., Berlin, Heidelberg, New York~: 
Springer-Verlag, 1983.

\m

\noindent{\sc Olver} (Peter J.)

[1986]
{\it Aplications of Lie Groups to Differential Equations}, Graduate Texts in Mathematics,
vol. 107, New York: Springer-Verlag, 1986; 2nd ed., revised, 1993.

[1986a] ``Noether's theorems and systems of Cauchy-Kovalevskaya type'', in {\it Nonlinear Systems
of Partial Differential Equations in Applied Mathematics}, Basil Nicolaenko, Darryl D. Holm, James M. Hyman, eds., Lectures
in Applied Mathematics, vol. 23, Part 2, Providence, RI: American Mathematical Society, 1986,
p. 81--104.

[2018] ``Emmy Noether's enduring legacy in symmetry'', {\it Symmetry: Culture and Science}, 29 (2018), p. 475--485.

\m

\noindent{\sc O'Raifeartaigh} (Lochlainn)

[1997] {\it The Dawning of Gauge Theory}, Princeton: Princeton University Press, 1997.

\m
\noindent{\sc Pauli} (Wolfgang)

[1921] \selectlanguage{german}``Relativit\"atstheorie'', {\it Encyklop\"adie der mathematischen Wissenschaften}\selectlanguage{english}, vol. V.2,
Leipzig: B. G. Teubner, 1921, p. 539--775; reprinted as a monograph, with a preface by Arnold
Sommerfeld, Leipzig: Teubner, 1921. {\it Collected Scientific Papers}, vol. 1, p. 1--237.


\newpage

\noindent{\sc Rowe} (David E.)

[1999] ``The G\"ottingen response to general relativity and Emmy Noether's theorems'', in 
{\it The Symbolic Universe, Geometry and Physics 1890--1930} (Milton Keynes, 1996), Jeremy J. Gray, ed.,  
Oxford, New York: Oxford University Press, 1999, p. 189--233.

\m

\noindent{\sc Souriau} (Jean-Marie)

[1964] {\it G\'eom\'etrie et relativit\'e}, Paris: Hermann, 1964.

[1970] {\it Structure des syst\`emes dynamiques}, Paris: Dunod, 1970; English
translation by C. H. Cushman-de Vries, {\it Structure of Dynamical Systems. A Symplectic
View of Physics}, Springer, 1997.

\m

\noindent{\sc Srinivasan} (Bhama) and {\sc Sally} (Judith D.), eds.

[1983] {\it Emmy Noether in Bryn Mawr}, New York, Berlin, Tokyo: Springer-Verlag, 1983.

\m

\noindent{\sc Stasheff} (James D.)

[1997] ``Deformation theory and the Batalin-Vilkovisky master equation'', in {\it Deformation Theory and Symplectic Geometry} (Ascona, 1996), Daniel Sternheimer, John Rawnsley, Simone Gutt, eds., Math. Phys. Stud., 20, Dordrecht: Kluwer Acad. Publ., 1997, p.~271--284.

\m

\noindent{\sc Trautman} (Andrzej)

[1967] ``Noether equations and conservation laws'', {\it Communications in Mathematical Physics}, 6
(1967), p. 248--261.

[1972] ``Invariance of Lagrangian Systems'', in {\it General Relativity}, papers in honor of J. L. Synge,
L. O'Raifeartaigh, ed., Oxford: Clarendon Press, 1972, p.~85--99.

\m

\noindent {\sc Utiyama} (Ryoyu )

[1947] ``On the interaction of mesons with the gravitational field. I.'', {\it Progress of Theoretical Physics}, vol. II, no. 2 (1947), p. 38--62.

\m

\noindent {\sc Vinogradov} (Alexandre M.)

[1977] ``On the algebro-geometric foundations of Lagrangian field theory'', {\it Doklady
Akad. Nauk SSSR}, 236 (1977), p. 284--287; English translation, {\it Soviet Mathematics Doklady}, 18
(1977), p. 120--1204.

\m
\noindent {\sc Vizgin} (Vladimir P.)

[1972] {\it Razvitije vzaimosvyazi principov invariantnosti s zakonami sokhranenija v klassicheskoj
fizike} [Russian] (Development of the Interconnection between Invariance Principles and Conservation
Laws in Classical Physics), Moscow: Nauka, 1972.

\m


\noindent{\sc Weyl} (Hermann)

[1918] {\it Raum, Zeit, Materie}, Berlin: J. Springer, 1918; 2nd ed., 1919; 3rd ed., revised, 1919,
and later editions.

[1935] ``Emmy Noether'' (Memorial address delivered in Goodhart Hall, Bryn Mawr College, on
April 26, 1935), {\it Scripta Mathematica}, 3 (1935), p. 201--220; {\it Gesammelte Abhandlungen}, vol. 3,
p. 425--444.
\m

\noindent{\sc Yang} (Chen Ning) and {\sc Mills} (Robert L.)

[1954] ``Conservation of isotopic spin and isotopic gauge invariance'', {\it Physical Review}, 96 (1954), p. 191--195.

\m

\noindent{\sc Zuckerman} (Gregg J.)

[1987] ``Action principles and global geometry'', in {\it Mathematical Aspects of String Theory} (University
of California, San Diego, 1986), S. T. Yau, ed., Advanced Series in Mathematical Physics,
vol. 1, Singapore, New Jersey, Hong Kong: World Scientific, 1987, p. 259--284.

\end{document}